\documentclass [11pt ,a4paper ,twoside]{article}
\usepackage{geometry}
\usepackage{amsmath}
\usepackage[english]{babel}
\usepackage[utf8]{inputenc}
\usepackage{fontenc}
\usepackage{newlfont}
\usepackage{indentfirst}
\usepackage{fancyhdr}
\usepackage{textcomp}
\usepackage{graphicx}
\usepackage{amsmath}
\usepackage{amssymb}
\usepackage{amsthm}
\usepackage{mathrsfs}
\usepackage{graphicx}
\usepackage{subfig}
\usepackage[dvips]{hyperref}
\theoremstyle{plain}
\newtheorem{Teo}{Theorem}[section]

\newtheorem{Prop}[Teo]{Proposition}
\theoremstyle{definition}
\newtheorem{Def}{Definition}
\newtheorem{Not}{Notation}
\theoremstyle{remark}
\newtheorem{Rem}{Remark}

\numberwithin{Teo}{section}
\numberwithin{Def}{section}
\numberwithin{Not}{section}
\numberwithin{Rem}{section}
\title{A note on a stochastic approach to Caffarelli-Silvestre Theorem}
\author{Cavina Michelangelo}
\begin{document}
\maketitle
\tableofcontents
\section*{Introduction}
\addcontentsline{toc}{section}{Introduction}
In this note we analyze the Caffarelli-Silvestre extension function using tools from the theory of stochastic analysis applied to Dirichlet problems. We use a stochastic approach to give the explicit formulation of the kernel associated to the Dirichlet problem which defines the Cafferelli-Silvestre extension function.\\

The connection between the Caffarelli-Silvestre extension and trace processes of diffusions in the upper half plane is known, and generally attributed to Molchanov and Ostrovskii (see \cite{MO}) in a more general context. Our aim here is giving a detailed and self contained proof of such results, which we could not find in literature.\\

Caffarelli and Silvestre proved in \cite{CS} that it is possible to represent the fractional Laplacian $(-\Delta)^s u$, for $0<s<1$, of a function $u \in C^2(\mathbb R^n)$ in terms of the solution $U \in C(\mathbb R_+^{n+1})$ to a (local) PDE problem in $\mathbb R^{n+1}_+:=\mathbb R^n \times (0,+\infty)$. It is possible to give an interpretation of this result based on the theory of stochastic analysis. Molchanov and Ostrovskii in \cite{MO} proved a probabilistic analogue of the extension technique, where they considered the trace process of a $2$-dimensional process $Z=(X,Y)$, where $X$ is a 1-dimensional Brownian motion, and $Y$ is a Bessel process, but they did not show the connection between the generator of the trace process and the boundary condition of the solution to a PDE problem. In this work we compute the value of the stochastic $s$-harmonic extension $U$ in $\mathbb R^{n+1}_+$ and we show that it is equal to the convolution between the boundary data $u$ and the expected Poisson-type kernel.\\

After proving the result and presenting it in two seminars we found a recent thesis about a generalization of the extension method used by Caffarelli and Silvestre. In his PhD thesis \cite{H} Herman showed that it is possible to generalize the extension method used in \cite{CS} to a wider family of non-local operators, using stochastic analysis and semigroup theory to prove that it is possible to represent a wide family of non-local operators in terms of the solution to a local PDE problem. The method used in Herman's work consists of considering the trace process of a of proper diffusion process in $\mathbb R^n \times[0,+\infty)$ and deriving the Neumann boundary conditions of a solution to a PDE from the generator of the trace process. The connection between the generator of the trace of a diffusion process and the Neumann boundary conditions was made in a stochastic sense by Hsu, see \cite{Hsu}. Roughly speaking, the connection is made by combining It\^o's formula and a random time change given by the inverse local time at the boundary. The Caffarelli-Silvestre extension technique can also be generalized to operators of the family $\{\varphi(-\Delta)\}$, where $\varphi$ is a complete Bernstein function and $-\Delta$ is the positive Laplace operator, see \cite{KM}.\\

Our approach is slightly simpler. Given a function $u \in C^2(\mathbb R^n)$ and point $(x_0,y_0) \in \mathbb R^{n+1}_+$ we consider a stochastic process $Z=(X,Y)$ starting from $(x_0,y_0)$, where $Y=\{Y_t\}_{t\geq0}$ is a Bessel process and $X=\{X_t\}_{t\geq 0}$ is a $n$-dimensional Brownian motion independent from $Y$, and we compute the expected value $\mathbb E^{(x_0,y_0)}\left[u\left(Z_{\tau_{\mathbb R^{n+1}_+}}\right)\right]=:w(x_0,y_0)$, here $\tau_{\mathbb R^{n+1}_+}$ denotes the first exit time for $Z$ from the domain $\mathbb R^{n+1}_+$. The function $w$ is the stochastic $s$-harmonic extension of $u$, and by the theorem about the stochastic solution to the Dirichlet problem (see \cite[theorem 9.2.14]{Oksendal}) the function $w$ satisfies the Dirichlet problem
\begin{equation}
\left(\text{D.P.}'\right) \begin{cases}
 \bigg(\frac{1-2s}{2y}\frac{\partial}{\partial y}+\frac 1 2 \Delta_{x,y}\bigg)w=0 \quad \text{in } R_+^{n+1},\\
w(x,0)=u(x) \quad\quad\quad\quad\quad\;\,\,\text{for }x \in \mathbb R^n,
\end{cases}
\end{equation}
which is the one associated to the Caffarelli-Silvestre extension function in \cite{CS}. Then, we prove that the function $w$ can be written under the form
\begin{equation}
w(x_0,y_0)=K_{y_0} * u (x_0)=\int_{\mathbb R^n}K_{y_0}(x_0-x)u(x)dx,
\end{equation}
where $K_y$ is the Poisson-type kernel
\begin{equation}\label{thesis theorem}
K_y(x)= \frac{1}{\pi^{\frac n 2}} \frac{\Gamma(s+\frac n 2)}{\Gamma(s)} \cdot \frac{y^{2s}}{(|x|^2 +y^2)^{\frac n 2 +s}}.
\end{equation}

For the notations and the theorems about the theory of stochastic analysis we reference \cite{Oksendal}. For the definition and properties of the Bessel process we reference \cite{Survey} and \cite{MO}.

\section{Stochastic Dirichlet problem}
In this section we list the definitions and theorems used in the proof of the Poisson-type Kernel formula given in section 3.
\begin{Not}
We will denote by $(\Omega, \mathcal F,P,\{\mathcal F_t\}_{t\geq 0})$ a filtered probability space $\Omega$ of $\sigma$-algebra $\mathcal F$, probability measure $P$ and filtration $\{\mathcal F_t\}_{t\geq 0}$. We will omit writing the filtration $\{\mathcal F_t\}_{t\geq 0}$ because we will always use the natural filtration associated to the Brownian motions mentioned in the following calculations.\\
\end{Not}

We reference \cite[chapters 7 and 9]{Oksendal} for the following definitions and theorems about diffusion processes.
\begin{Def}[It\^o diffusion {\cite[definition 7.1.1]{Oksendal}}]
A (time homogeneous) It\^o diffusion is a stochastic process
\begin{align}
X:[0,+\infty) \times \Omega & \longrightarrow \mathbb R^n \\
(t,\omega) & \longmapsto X_t(\omega) \nonumber
\end{align}
satisfying the following stochastic differential equation:
\begin{equation}
dX_t = b(X_t)dt + \sigma(X_t) dB_t; \quad X_0=x.
\end{equation}
Here $x \in \mathbb R^n$ is the starting point at the time $t=0$, $B=\{B_t\}_{t\geq0}$ is a standard $m$-dimensional Brownian motion, and
\begin{equation}
b:\mathbb R^n \longrightarrow \mathbb R^n; \quad \sigma: \mathbb R^n \longrightarrow \mathbb R^{n \times m}
\end{equation}
are coefficients satisfying proper conditions (see \cite{Oksendal}, chapter 7).\\
Let $f:\mathbb R^n \longrightarrow \mathbb R$. We denote by $\mathbb E^x\left[f(X_t)\right]$ the expected value (w.r.t. the probability measure P) of the function $f$ evaluated at the It\^o diffusion $X$ of starting point $x \in \mathbb R^n$ at the time $t\geq0$.
\end{Def}
\begin{Def}[Infinitesimal generator {\cite[definition 7.3.1]{Oksendal}}]
Let $X=\{X_t\}_{t\geq0}$ be an It\^o diffusion in $\mathbb R^n$. The infinitesimal generator $\mathcal A$ of $X$ is defined by
\begin{equation}
\mathcal A f(x)= \lim_{t\rightarrow 0^+} \frac{\mathbb E^x\left[f(X_t)\right]-f(x)}{t}, \quad \text{for } x \in \mathbb R^n.
\end{equation}
\end{Def}
The operator $\mathcal A$ is well defined everywhere for all the functions $f \in C_0^2(\mathbb R^n)$.
\begin{Teo}[Characterization of infinitesimal generators {\cite[theorem 7.3.3]{Oksendal}}]\label{Characterization of infinitesimal generators theorem}
Let $\{X_t\}_{t\geq0}$ be an It\^o diffusion satisfying
\begin{equation*}
dX_t = b(X_t)dt + \sigma(X_t) dB_t.
\end{equation*}
Let $f \in C_0^2(\mathbb R^n)$. Then
\begin{equation}
\mathcal A f(x) = \sum_{i=1}^n b_i(x)\frac{\partial}{\partial x_i}f(x) + \frac 1 2 \sum_{i,j=1}^n \left( \sigma \cdot \sigma^T \right)_{i,j}(x) \frac{\partial^2}{\partial x_i \partial x_j}f(x).
\end{equation}
Here $\left( \sigma \cdot \sigma^T \right)_{i,j}$ denotes the component of coordinates $(i,j)$ of the matrix $ \sigma \cdot \sigma^T$, where $\sigma^T$ is the transposed of $\sigma$.
\end{Teo}
\begin{Def}[First exit time for a stochastic process]
Let $D \subseteq \mathbb R^n$, let $X:[0,+\infty) \times \Omega \longmapsto \mathbb R^n$ be an It\^o diffusion. We denote by first exit time of $X$ from $D$ the random variable
\begin{equation}
\tau_D: \Omega \longmapsto [0,+\infty]; \quad \tau_D(\omega):= \inf\{t>0 \; | \; X_t(\omega) \not \in D\}.
\end{equation}
Moreover, we denote by $X$ at the time $\tau_D$ the random variable
\begin{equation}
X_{\tau_D} : \Omega \longrightarrow \mathbb R^n; \quad X_{\tau_D}(\omega):=\begin{cases}
X_{\tau_D(\omega)}(\omega) \; \text{ if } \tau_D(\omega)<+\infty,\\
0 \quad \quad\quad \quad\;\;\text{otherwise.}
\end{cases}
\end{equation}
\end{Def}
\begin{Def}[Regularity with respect to an It\^o diffusion {\cite[definition 9.2.8]{Oksendal}}]
Under the previous notations, assume $X$ is an It\^o diffusion. We say that a point $y \in \partial D$ is regular w.r.t. $X$ if
\begin{equation}
P^y[\tau_D=0]=1,
\end{equation}
otherwise $y$ is called irregular. Here $P^y[\tau_D=0]$ denotes the probability that the first exit time of the diffusion $X$ starting at the point $y$ is equal to 0.
\end{Def}
\begin{Teo}[Stochastic solution to the Dirichlet problem {\cite[theorem 9.2.14]{Oksendal}}]\label{Stochastic solution of the Dirichlet problem theorem}
Let $D \subseteq \mathbb R^n$, let $u \in C(\partial D)$, $u$ bounded.
Consider the Dirichlet problem
\begin{equation}
\left( \text{D.P.} \right) \begin{cases}
\mathcal A w=0 \quad\quad\quad \text{in } \partial D,\\
w(x)=u(x) \quad\,\text{for }x \in \partial D.
\end{cases}
\end{equation}
Let $\{X_t\}_{t\geq 0}$ be an It\^o diffusion such that the infinitesimal generator of $\{X_t\}_{t\geq 0}$ is $\mathcal A$.\\
Consider the function
\begin{align}
f:\overline{D} \longrightarrow \mathbb R^n\nonumber \\
f(x)=\mathbb E^x\left[u(X_{\tau_D})\right].
\end{align}
Then, under suitable hypotheses, $f$ is a solution to
\begin{equation}
\left(\text{D.P.}'\right) \begin{cases}
\mathcal A f=0 \quad\quad\quad \text{in } D,\\
f(x)=u(x) \quad\, \text{for } x \in \partial D, \; \; x  \text{ regular w.r.t. } \{X_t\}_{t\geq 0}.
\end{cases}
\end{equation}
\end{Teo}
\section{The Bessel process}
In this section we define the Bessel process and enunciate the properties we use in the proof of the Poisson-type Kernel formula given in section 3.
\begin{Def}
Let $0<s<1$ We denote by Bessel process in $\mathbb R$ the stochastic process
\begin{align*}
Y:[0,+\infty)\times \Omega_2 &\longrightarrow \mathbb R\\
(t,\omega)&\longmapsto Y_t(\omega),
\end{align*}
satisfying the following stochastic differential equation:
\begin{equation}
dY_t= \frac{1-2s}{2Y_t}dt + dB^{n+1}_t,
\end{equation}
where $B^{n+1}_t$ is a 1-dimensional Brownian motion.
\end{Def}
\begin{Prop}
The following facts about the Bessel process $\{Y_t\}_{t\geq 0}$ hold (see \cite{Survey}, \cite{MO}):
\begin{enumerate}
\item  $\{Y_t\}_{t\geq 0}$ is a continuous diffusion process.
\item Let $y_0$ be a starting point. Then the trajectories $t \mapsto Y_t(\omega)$ hit the point 0 in a finite amount of time almost surely.
\item  Let $y_0$ be a starting point. The random variable ``first hitting time for the process $\{Y_t\}_{t\geq 0}$ starting from $y_0$ and hitting 0" has a density with respect to the Lebesque measure (see \cite{Survey}, page 8, equation (15)). The density function is
\begin{equation}\label{density of Bessel hitting time}
\varPhi_{y_0}(t) = \chi_{(0,+\infty)}(t) \frac{1}{t \Gamma(s)}\left(\frac{{y_0}^2}{2t}\right)^s e^{-\frac{{y_0}^2}{2t}}.
\end{equation}
\end{enumerate}
\end{Prop}
\begin{Rem}
For any choice of $\alpha >0$ and $M>0$ we have
\begin{equation}\label{inner integral is equal to 1}
\int_0^{+\infty} \frac{1}{t \Gamma(\alpha)}\left(\frac{M}{2t}\right)^{\alpha} e^{-\frac{M}{2t}}dt=1.
\end{equation}
\end{Rem}
\section{Caffarelli-Silvestre theorem}
In this section we prove that the Poisson Kernel formula associated to Caffarelli-Silvestre theorem can be obtained using theorem \ref{Stochastic solution of the Dirichlet problem theorem}.\\
We begin by recalling the definition of fractional Laplacian and Caffarelli-Silvestre theorem.
\begin{Def}[Fractional Laplacian]
Let $0<s<1$. Let $u \in C^2(\mathbb R^n)$, u bounded. We define the fractional Laplacian
\begin{equation}
(-\Delta)^s u(x_0):=A_{n,s} \cdot \text{P.V.} \int_{\mathbb R^n} \frac{u(x_0)-u(x)}{|x_0-x| ^{n+2s}} dx.
\end{equation}
Here $A_{n,s}$ is a constant depending only on $n$ and $s$.
\end{Def}
\begin{Teo}[Caffarelli-Silvestre]
Let $D= \mathbb R_+^{n+1}:= \mathbb R^n \times (0,+\infty) \ni (x,y)$. We identify $\partial \mathbb R_+^{n+1} \equiv \mathbb R^n$.\\
Let $u \in C^2(\mathbb R^n)$, $u$ bounded. Let $U:\overline{D} \rightarrow \mathbb R$ be the solution to
\begin{equation}
\left( \text{D.P.} \right) \begin{cases}
 \text{div}\left(y^{1-2s}\nabla U\right)=0 \quad \text{in } D =\mathbb R_+^{n+1},\\
U(x,0)=u(x) \quad\quad\;\;\,\,\text{for }x \in \mathbb R^n,
\end{cases}
\end{equation}
which is equivalent to
\begin{equation}\label{DP' in caffarelli silvestre}
\left(\text{D.P.}'\right) \begin{cases}
 \bigg(\frac{1-2s}{2y}\frac{\partial}{\partial y}+\frac 1 2 \Delta_{x,y}\bigg)U=0 \quad \text{in } D =\mathbb R_+^{n+1},\\
U(x,0)=u(x) \quad\quad\quad\quad\quad\;\,\,\text{for }x \in \mathbb R^n.
\end{cases}
\end{equation}
Then
\begin{equation}
(-\Delta)^s u(x)= - A_{n,s} \lim_{y \rightarrow 0^+} y^{1-2s} \cdot  \frac{\partial}{\partial_y }U(x,y).
\end{equation}
\end{Teo}
Caffarelli and Silvestre proved (see \cite{CS}, section 3) that this theorem follows from the following formula about a Poisson-type kernel.
\begin{Prop}[Poisson-type Kernel formula]\label{Poisson-type Kernel formula proposition}
Let $u \in C^2(\mathbb R^n)$, $u$ bounded. Let $D=\mathbb R_+^{n+1}$, $D \equiv\mathbb R^n$.\\
Consider the following Dirichlet problem
\begin{equation}\label{Dirichlet problem in theorem}
\left(\text{D.P.}'\right) \begin{cases}
 \bigg(\frac{1-2s}{2y}\frac{\partial}{\partial y}+\frac 1 2 \Delta_{x,y}\bigg)\phi=0 \quad \text{in } D =\mathbb R_+^{n+1},\\
\phi(x,0)=u(x) \quad\quad\quad\quad\quad\;\,\,\text{for }x \in \mathbb R^n.
\end{cases}
\end{equation}
Then the function
\begin{equation}
U(x_0,y_0)=K_{y_0} * u (x_0)=\int_{\mathbb R^n}K_{y_0}(x_0-x)u(x)dx,
\end{equation}
is a solution to (\ref{Dirichlet problem in theorem}), where
\begin{equation}
K_y(x)=C_{n,s} \cdot \frac{y^{2s}}{(|x|^2 +y^2)^{\frac n 2 +s}}.
\end{equation}
Here the constant $C_{n,s}$ is
\begin{equation}
C_{n,s}:= \frac{1}{\pi^{\frac n 2}} \frac{\Gamma(s+\frac n 2)}{\Gamma(s)},
\end{equation}
\end{Prop}

We are going to give a proof of proposition \ref{Poisson-type Kernel formula proposition} using the results from section 1.
\proof
Consider the It\^o diffusion $Z=\{Z_t\}_{t\geq0}$ in $\mathbb R^{n+1}$ satisfying
\begin{equation}
\left(\begin{array}{c} dZ_t^1\\dZ_t^2\\ \vdots \\dZ_t^n\\dZ_t^{n+1}  \end{array}\right)=:\left(\begin{array}{c} dX_t^1\\dX_t^2\\ \vdots \\dX_t^n\\dY_t  \end{array}\right)=\left(\begin{array}{c}0\\0\\ \vdots \\0\\ \frac{1-2s}{2Y_t}  \end{array}\right)dt +I_{n+1}\cdot \left(\begin{array}{c} dB_t^1\\dB_t^2\\ \vdots \\dB_t^n\\dB_t^{n+1}  \end{array}\right),
\end{equation}
where $(B_t^1,\dots,B_t^{n+1})$ is a standard $(n+1)$-dimensional Brownian motion, and $I_{n+1} \in \mathbb R^{(n+1) \times(n+1)}$ denotes the identity matrix
\begin{equation}
I_{n+1}:=\left(\begin{array}{ccccc} 1&0&0&\dots &0 \\ 0&1 &0& \dots&0\\ 0&0&1&\dots &0\\\vdots&\vdots&\vdots&\ddots&0\\0&0&0&\dots &1 \end{array}\right).
\end{equation}
Let $z_0=(x_0,y_0)$ be the starting point of $Z$. The process $Z$ can be written as $Z=(X,Y)$, where:
\begin{itemize}
\item $X=\{X_t\}_{t\geq0}$ is a standard $n$-dimensional Brownian motion starting from $x_0$.
\item $Y=\{Y_t\}_{t\geq0}$ is a Bessel process starting from $y_0$ and independent from $X$.
\end{itemize}
Using theorem \ref{Characterization of infinitesimal generators theorem} we get that the infinitesimal generator of the process $Z$ is
\begin{equation}
\mathcal A f(x,y)=\left( \frac{1-2s}{2y} \frac{\partial}{\partial y} + \frac 1 2 \Delta_{x,y}\right) f(x,y).
\end{equation}
$\mathcal A$ is the operator associated to (\ref{DP' in caffarelli silvestre}).\\
Moreover, all the points $z=(x,0) \in \partial \mathbb R^{n+1}_+$ are regular with respect to $\{Z_t\}_{t\geq 0}$ because, when $0<s<1$, the Bessel process $Y$ oscillates around 0 and hits it infinitely many times, in every interval of time starting from the time of first hitting 0, with probability 1 (see \cite{MO}, \cite{Survey}).\\
So we may apply theorem \ref{Stochastic solution of the Dirichlet problem theorem} and get that the function
\begin{equation}
w(x_0,y_0)=\mathbb E^{x_0,y_0}\left[ u(Z_{\tau_{D}})\right] \quad \text{for } (x_0,y_0) \in \overline{\mathbb R^{n+1}_+}
\end{equation}
satisfies
\begin{equation}
\left(\text{D.P.}'\right) \begin{cases}
 \bigg(\frac{1-2s}{2y}\frac{\partial}{\partial y}+\frac 1 2 \Delta_{x,y}\bigg)w=0 \quad \text{in } D =\mathbb R_+^{n+1},\\
w(x,0)=u(x) \quad\quad\quad\quad\quad\;\,\,\text{for }x \in \mathbb R^n.
\end{cases}
\end{equation}
Now we compute $w$ for a starting point $(x_0,y_0)$, $y_0 >0$.\\
Let the probability space $(\Omega_1,{\mathcal F}_1, P_1)$ be the domain of the random variables $X_t:\Omega_1 \rightarrow \mathbb R^n$, and let the probability space $(\Omega_2,{\mathcal F}_2, P_2)$ be the domain of the random variables $Y_t:\Omega_2 \rightarrow \mathbb R^n$.\\
Then the domain of the random variables $Z_t=(X_t,Y_t)$ is the product space $(\Omega_1 \times \Omega_2,{\mathcal F}_1 \otimes {\mathcal F}_2, P_1 \times P_2)$. We are going to compute
\begin{equation}
w(x_0,y_0)=\int_{\Omega_1 \times \Omega_2}u\left(Z^{(x_0,y_0)}_{\tau_{D}(\omega_1,\omega_2)}(\omega_1,\omega_2)\right)d(P_1 \times P_2)(\omega_1,\omega_2).
\end{equation}
Now we observe that $Y_{\tau_{D}}=0$ almost surely, because the process $Y$ is continuous and $\partial D=\{(x,y) \in \mathbb R^n \times \mathbb R\; | \; y=0 \}$, $X$ is a Brownian motion independent from $\omega_2$, and the exit time $\tau_{D}$ doesn't depend on $\omega_1$ because the process $(X,Y)$ exits from the domain $D$ if and only if the process $Y$ exits from $(0,+\infty)$.\\
So, with a little abuse of notation, we may write
\begin{equation}
u\left(Z^{(x_0,y_0)}_{\tau_{D}(\omega_1,\omega_2)}(\omega_1,\omega_2)\right)=u\left(X^{x_0}_{\tau_{D}(\omega_2)}(\omega_1)\right).
\end{equation}
We apply Fubini-Tonelli theorem and we get
\begin{equation}
w(x_0,y_0)=\int_{\Omega_2}\left[\int_{\Omega_1}u\left(X^{x_0}_{\tau_{D}(\omega_2)}(\omega_1)\right) dP_1(\omega_1)\right] dP_2(\omega_2).
\end{equation}
However, $X$ is a Brownian motion, so 
\begin{equation}
X^{x_0}_{\tau_{D}(\omega_2)}\sim \mathcal N(x_0, \tau_{D}(\omega_2) \cdot I_n),
\end{equation}
i.e. $X$ has the same probability distribution as a multivariate normal variable of mean value equal to the vector $x_0$, and matrix of covariances equal to $ \tau_{D}(\omega_2) \cdot I_n$.\\
So we use the equation of the density of the multivariate normal variable to get
\begin{equation}
w(x_0,y_0)=\int_{\Omega_2}\left[\int_{\mathbb R^n}\frac{1}{(2 \pi \tau_{D}(\omega_2))^{\frac n 2}}e^{-\frac{|x-x_0|^2}{2\tau_{D}(\Omega_2)}}u(x) dx\right] dP_2(\omega_2).
\end{equation}
Now we use the density of the variable $\tau_{D}$ from equation (\ref{density of Bessel hitting time}) to get
\begin{equation}
w(x_0,y_0)=\int_0^{+\infty}  \frac{1}{t \Gamma(s)}\left(\frac{{y_0}^2}{2t}\right)^s e^{-\frac{{y_0}^2}{2t}}\left[\int_{\mathbb R^n}\frac{1}{(2 \pi t)^{\frac n 2}}e^{-\frac{|x-x_0|^2}{2t}}u(x) dx\right] dt.
\end{equation}
We change the order of integration and rearrange the factors to get
\begin{align}
w(x_0,y_0)=&\int_{\mathbb R^n}  \frac{1}{\pi^{\frac n 2}}\frac{\Gamma(s+\frac n 2)}{\Gamma(s)}\frac{{y_0}^{2s}}{(|x-x_0|^2+{y_0}^2)^{s+\frac n 2}}u(x) \;\cdot \\
&\left[\int_0^{+\infty}\frac{1}{t \Gamma(s+\frac n 2)}\left(\frac{|x-x_0|^2+{y_0}^2}{2t}\right)^{s+\frac n 2}e^{-\frac{|x-x_0|^2+{y_0}^2}{2t}}dt\right] dx.\nonumber
\end{align}
However, by equation (\ref{inner integral is equal to 1}) with $M=|x-x_0|^2+{y_0}^2$ and $\alpha=s+n/2$, we get
\begin{equation}
w(x_0,y_0)=\int_{\mathbb R^n}  \frac{1}{\pi^{\frac n 2}} \frac{\Gamma(s+\frac n 2)}{\Gamma(s)}\frac{{y_0}^{2s}}{(|x_0-x|^2+{y_0}^2)^{s+\frac n 2}}u(x)dx.
\end{equation}
We define
\begin{equation}
K_y(x):= C_{n,s}\cdot \frac{{y}^{2s}}{(|x|^2+{y}^2)^{s+\frac n 2}},
\end{equation}
where
\begin{equation}
C_{n,s}:= \frac{1}{\pi^{\frac n 2}} \frac{\Gamma(s+\frac n 2)}{\Gamma(s)},
\end{equation}
and we get
\begin{equation}
w(x_0,y_0)=\int_{\mathbb R^n} K_{y_0}(x_0-x)u(x)dx=K_{y_0} * u (x_0)=U(x_0,y_0).
\end{equation}
So we proved that $w\equiv U$, and that $w$ is a solution to the Dirichlet problem (\ref{Dirichlet problem in theorem}), so the statement is proved.
\endproof


\begin{thebibliography}{7}
\addcontentsline{toc}{section}{References}
\bibitem{CS} L. Caffarelli; L. Silvestre. \textit{An extension problem related to the fractional Laplacian.} Communications in partial differential equations 32, 8 (2007), 1245$-$1260.
\bibitem{Survey}  A. G\"oing-Jaeschke; M. Yor. \textit{A survey and some generalizations of Bessel processes.} Bernoulli 9, no. 2 (2003), 313$-$349.
\bibitem{H} J.A. Herman. \textit{The Harmonic Extension Technique with Applications to Optimal Stopping.} PhD thesis, University of Warwick (2020).
\bibitem{Hsu} P. Hsu. \textit{On excursions of reflecting brownian motion.} Transactions of the American Mathematical Society 296, 1 (1986), 239$-$264.
\bibitem{KM} M. Kwa\'snicki; J. Mucha. \textit{Extension technique for complete Bernstein functions of the Laplace operator.} Journal of Evolution Equations (2017), 1$-$39.
\bibitem{MO} S.A. Molchanov; E. Ostrovskii. \textit{Symmetric stable processes as traces of degenerate diffusion processes.} Theory of Probability \& Its Applications 14, 1 (1969), 128$-$131.
\bibitem{Oksendal} B. \O ksendal. \textit{Stochastic differential equations.} An introduction with applications. Fifth edition. Universitext. Springer-Verlag, Berlin (1998).
\end{thebibliography}
\end{document}